\documentclass[10pt]{article}

\usepackage{amsmath, amssymb, amscd, fancyhdr}

\newcommand{\version}{version 1.0,\ \ Oct. 24, 2012}

\renewcommand{\leftmark}{{\scriptsize Rational curves and special metrics on twistor spaces}}

\def\eqref#1{(\ref{#1})}
\newcommand{\goth}{\mathfrak}

\newcommand{\arrow}{{\:\longrightarrow\:}}

\newcommand{\C}{{\Bbb C}}
\newcommand{\R}{{\Bbb R}}

\def\1{\sqrt{-1}\:}
\newcommand{\restrict}[1]{{\left|_{{\phantom{|}\!\!}_{#1}}\right.}}
\newcommand{\cntrct}                
{\hspace{2pt}\raisebox{1pt}{\text{$\lrcorner$}}\hspace{2pt}}

\def\Bbb#1{\mathbb #1}

\newcommand{\calo}{{\cal O}}

\renewcommand{\bar}{\overline}
\renewcommand{\phi}{\varphi}
\renewcommand{\epsilon}{\varepsilon}
\renewcommand{\geq}{\geqslant}

\newcommand{\im}{\operatorname{im}}

\newcommand{\Tot}{\operatorname{Tot}}

\newcommand{\Vol}{\operatorname{Vol}}

\newcommand{\Lie}{\operatorname{Lie}}
\newcommand{\Tw}{\operatorname{Tw}}

\newcommand{\comment}[1]{{}}

\def\blacksquare{\hbox{\vrule width 4pt height 4pt depth 0pt}}
\def\endproof{\blacksquare}

\makeatletter

\newcounter{Mycounter}[section]
\newcounter{lemma}[section]
\renewcommand{\thelemma}{\noindent{Lemma \thesection.\arabic{lemma}}}
\newcommand{\lemma}{%
     \setcounter{lemma}{\value{Mycounter}}
     \refstepcounter{lemma}
     \stepcounter{Mycounter}
     {\bf \thelemma:\ }}

\newcounter{claim}[section]
\renewcommand{\theclaim}{\noindent{Claim \thesection.\arabic{claim}}}
\newcommand{\claim}{%
     \setcounter{claim}{\value{Mycounter}}
     \refstepcounter{claim}
     \stepcounter{Mycounter}
     {\bf \theclaim:\ }}

\newcounter{sublemma}[section]
\newcounter{corollary}[section]
\renewcommand{\thecorollary}{\noindent{Corollary \thesection.\arabic{corollary}}}
\newcommand{\corollary}{%
     \setcounter{corollary}{\value{Mycounter}}
     \refstepcounter{corollary}
     \stepcounter{Mycounter}
     {\bf \thecorollary:\ }}

\newcounter{theorem}[section]
\renewcommand{\thetheorem}{\noindent{Theorem \thesection.\arabic{theorem}}}
\newcommand{\theorem}{%
     \setcounter{theorem}{\value{Mycounter}}
     \refstepcounter{theorem}
     \stepcounter{Mycounter}
     {\bf \thetheorem:\ }}

\newcounter{conjecture}[section]
\newcounter{proposition}[section]
\renewcommand{\theproposition}
       {\noindent{Proposition \thesection.\arabic{proposition}}}
\newcommand{\proposition}{%
     \setcounter{proposition}{\value{Mycounter}}
     \refstepcounter{proposition}
     \stepcounter{Mycounter}
     {\bf \theproposition:\ }}

\newcounter{definition}[section]
\renewcommand{\thedefinition}
       {\noindent{Definition~\thesection.\arabic{definition}}}
\newcommand{\definition}{%
     \setcounter{definition}{\value{Mycounter}}
     \refstepcounter{definition}
     \stepcounter{Mycounter}
     {\bf \thedefinition:\ }}

\newcounter{example}[section]
\newcounter{remark}[section]
\renewcommand{\theremark}{\noindent{Remark \thesection.\arabic{remark}}}
\newcommand{\remark}{%
     \setcounter{remark}{\value{Mycounter}}
     \refstepcounter{remark}
     \stepcounter{Mycounter}
     {\bf \theremark:\ }}

\newcounter{problem}[section]
\newcounter{question}[section]
\@addtoreset{equation}{section}
\@addtoreset{footnote}{section}
\makeatother

\begin{document}

\begin{center}
{\LARGE\bf Rational curves and special metrics on twistor spaces
}
\\[4mm]
Misha Verbitsky\footnote{Partially supported by RFBR grants
 12-01-00944-Á,  10-01-93113-NCNIL-a, and
AG Laboratory NRI-HSE, RF government grant, ag. 11.G34.31.0023.}
\\[4mm]

{\tt verbit@verbit.ru}

\hfill

{\em \hfil To Professor H. Blaine Lawson at his 70th birthday}

\end{center}

{\small 
\hspace{0.15\linewidth}
\begin{minipage}[t]{0.7\linewidth}
{\bf Abstract} \\
A Hermitian metric $\omega$ on a complex manifold
is called {\bf SKT} or {\bf pluriclosed} if $dd^c\omega=0$.
 Let $M$ be a twistor space of 
a compact, anti-selfdual Riemannian manifold,
admitting a pluriclosed Hermitian metric. We prove that
in this case $M$ is K\"ahler, hence isomorphic
to $\C P^3$ or a flag space. This result is obtained from 
rational connectedness of the twistor space, due to F. Campana.
As an aside, we prove that the moduli space of
rational curves on the twistor space of a K3 surface
is Stein. 
\end{minipage}
}

{
\small
\tableofcontents
}


\section{Introduction}


\subsection{Special Hermitian metrics on complex manifolds}

The world of non-K\"ahler complex geometry is infinitely 
bigger than that inhabited by K\"ahler manifolds. For instance,
as shown by Taubes \cite{_Taubes_} (see also \cite{_Panov_Petrunin_}, 
any finitely-generated group can be realized
as a fundamental group of a compact complex 
manifold. On contrary, the K\"ahler condition
puts big restrictions on the fundamental group.

However, there are not many constructions
which lead to explicit non-K\"ahler complex
manifolds. There are many homogeneous and locally
homogeneous manifolds (such as complex nilmanifolds),
which are known to be non-K\"ahler. The locally
conformally K\"ahler manifolds are non-K\"ahler
by a theorem of Vaisman (\cite{_Vaisman_}).
Kodaira class VII surfaces (forming a vast and
still not completely understood class 
of complex surfaces) are never K\"ahler.
Finally, the twistor spaces, as shown by
Hitchin, are never K\"ahler, except two examples:
$\C P^3$, being a twistor space of $S^4$, and
the flag space, being a twistor space of $\C P^2$
(\cite{_Hitchin_}).

There are many ways to weaken the K\"ahler condition
$d\omega=0$.
Given a Hermitian form $\omega$ on a complex $n$-manifold, 
one may consider an equation $d(\omega^k)=0$. For
$1 < k < n-1$, this equation is equivalent to $d\omega=0$, but the equation
$d(\omega^{n-1})$ is quite non-trivial.
Such metrics are called {\bf balanced}.
All twistor spaces are balanced (\cite{_Michelsohn_});
also, all Moishezon manifolds are balanced
(\cite{_Alessandrini_Bassanellu:bimero_}).

Another way to weaken the K\"ahler condition
is to consider the equation $dd^c(\omega^k)=0$, where $d^c=-IdI$; this 
equation is non-trivial for all $0<k <n$.
When $k=1$, a metric satisfying $dd^c \omega=0$
is called {\bf pluriclosed}, or {\bf strong K\"ahler
torsion (SKT) metric}; such metrics are quite important in physics
and in generalized complex geometry.

A Hermitian metric satisfying  $dd^c(\omega^{n-1})=0$
is called {\bf Gauduchon}. As shown by P. Gauduchon
(\cite{_Gauduchon_1984_}), every Hermitian metric
is conformally equivalent to a Gauduchon metric,
which is unique in its conformal class up to a constant multiplier.

Since a twistor space has complex dimension 3, 
and is balanced, the only non-trivial metric
condition (among those mentioned above) 
for the twistor space is $dd^c(\omega)=0$. 

The main result of this paper is the following
theorem, similar to Hitchin's theorem on non-K\"ahlerianity
of twistor spaces.

\hfill

\theorem
Let $M$ be a twistor space of a compact
4-dimensional anti-selfdual Riemannian
manifold. Assume that $M$ admits a 
pluriclosed Hermitian form $\omega$: $dd^c(\omega)=0.$
Then $M$ is Kaehler.

{\bf Proof:} \ref{_plurica_then_Kah_Corollary_}. \endproof

\subsection{Strongly Gauduchon and symplectic Hermitian metrics}

The Gauduchon, pluriclosed and all the 
rest of the $dd^c(\omega^k)=0$ Hermitian metrics
have an interesting variation of a cohomological nature.

\hfill

\definition
Let $(M,I)$ be a complex manifold, and $\omega$
a Hermitian form.
We say that $\omega^k$ is {\bf strongly pluriclosed}
if any of the following equivalent conditions are
satisfied.
\begin{description}
\item[(i)] $d(\omega^k)$ is $dd^c$-exact.
\item[(ii)] $\omega^k$ is a $(k,k)$-part of a 
closed $2k$-form.
\end{description}
Notice that either of these conditions
easily implies $dd^c(\omega^k)=0$, but these
conditions are significantly stronger.

\hfill

For $k=1$ and $n-1$ this condition is especially interesting.
When a pluriclosed Hermitian form $\omega$ is (1,1)-part
of a closed (and hence symplectic) form $\tilde \omega$, $\omega$ is
called {\bf taming} or {\bf Hermitian symplectic},
and when $(\omega)^{n-1}$ is $(n-1,n-1)$-part of
a closed form, $\omega$ is called {\bf strongly Gauduchon}
(\cite{_Popovici_}).

\hfill

In the paper \cite{_Streets_Tian_}
Streets and Tian have constructed a 
parabolic flow for Hermitian symplectic metric,
analoguous to the K\"ahler-Ricci flow. They asked
whether there exists a compact complex
Hermitian symplectic manifold not admitting
a K\"ahler structure. This question was considered
in \cite{_Enrietti_Grancharov_Fino_} and \cite{_Enrietti_Fino_Vezzoni_}
for complex nilmanifolds. In \cite{_Enrietti_Fino_Vezzoni_}
it was shown that complex nilmanifolds cannot admit
Hermitian symplectic metrics.
However, the pluriclosed metrics exist 
on many complex nilmanifolds.

\hfill

The present paper grew as an attempt to answer the Streets-Tian's
question for twistor spaces. However, it was found that the
twistor spaces are not only never Hermitian symplectic,
they never admit a pluriclosed metric unless 
they are K\"ahler.

\subsection{Rational curves and pluriclosed metrics}

The results of the present paper are based on 
the study of the moduli of rational curves. Unlike
many complex non-algebraic manifolds, the twistor spaces
are very rich in curves: there exists a smooth rational
curve passing through any finite subset of a twistor
space (\ref{_rationally_conne_Claim_}).

For an almost complex structure $I$
equipped with a taming symplectic form,
all components of the space of complex curves 
are compact, by Gromov's compactness theorem
(\cite{_Gromov:curves_,_Audin_}).  I will show that
the same is true for pluriclosed metrics, if
$I$ is integrable (\ref{_compa_from_plurine_Corollary_}). 
This is used to prove that
a twistor space admitting a pluriclosed metric
is actually Moishezon (\ref{_Moishezon_if_quasilines_Theorem_}).

However, Moishezon varieties satisfy the $dd^c$-lemma.
This is used to show that any pluriclosed metric
is in fact Hermitian symplectic.

Finally, by using the Peternell's theorem from 
\cite{_Peternell_}, we prove that no Moishezon manifold
can be Hermitian symplectic (\ref{_plurica_then_Kah_Corollary_}).

\section{Twistor spaces for 4-dimensional Riemannian manifolds
and the space of rational curves}

\subsection{Twistor spaces for 4-dimensional Riemannian manifolds: 
definition and basic results}

\definition
Let $M$ be a Riemannian 4-manifold. Consider the action of
the Hodge $*$-operator: $*:\; \Lambda^2 M \arrow \Lambda^2 M$.
Since $*^2 =1$, the eigenvalues are $\pm 1$, and one
has a decomposition $\Lambda^2 M = \Lambda^+ M \oplus \Lambda^- M$
onto {\bf selfdual} ($*\eta=\eta$) and 
{\bf anti-selfdual} ($*\eta=-\eta$) forms.

\hfill

\remark If one changes the orientation of $M$,
leaving  metric the same, $\Lambda^+ M$ and $\Lambda^- M$
are exchanged. Therefore, $\dim \Lambda^2 M=6$ 
implies $\dim \Lambda^\pm(M)=3$.

\hfill

\remark
Using the isomorphism $\Lambda^2 M = \goth{so} (TM)$,
we interpret $\eta \in \Lambda^2_m M$ as an endomorphisms of $T_mM$.
Then the unit vectors $\eta \in \Lambda^+_mM$ correspond
to oriented, orthogonal complex structures on $T_m M$.

\hfill

\definition
Let $\Tw(M):=S\Lambda^+ M$ be the set of unit vectors in $\Lambda^+M$.
At each point $(m,s)\in \Tw(M)$, consider the decomposition
$T_{m,s}\Tw(M)= T_m M \oplus T_sS\Lambda^+_m M$, induced by the
Levi-Civita connection.
Let $I_s$ be the complex structure on $T_mM$ induced by $s$,
$I_{S\Lambda^+_m M}$ the complex structure on $S\Lambda^+_m M=S^2$
induced by the metrics and orientation, and 
${\cal I}:\; T_{m,s}\Tw(M)\arrow T_{m,s}\Tw(M)$
be equal to ${\cal I}_s \oplus I_{S\Lambda^+_m M}$.
An almost complex manifold $(\Tw(M), {\cal I})$
is called {\bf the twistor space} of $M$.

\hfill

The following results about twistor spaces are
well known (see e.g. \cite{_Besse:Einst_Manifo_}).

\hfill

\theorem
The almost complex structure on $(\Tw(M), {\cal I})$
is a conformal invariant of $M$. Moreover, one can
reconstruct the conformal structure on $M$ from the
almost complex structure on $\Tw(M)$ and its 
anticomplex involution $(m,s)\arrow (m, -s)$.
\endproof

\hfill

\theorem
$(\Tw(M), {\cal I})$ is a complex manifold if
  and only if $W^+=0$, where $W^+$ (``self-dual conformal
curvature'') 
is an autodual component of the curvature tensor. Such manifolds are
called {\bf conformally half-flat} or {\bf ASD
(anti-selfdual).} \endproof

\subsection{Rational curves on $\Tw(M)$}

\definition
{\bf An ample rational curve} on a complex manifold $M$ is
a smooth curve $S \cong \C P^1\subset M$ such that 
$NS=\bigoplus_{k=1}^{n-1}\calo(i_k)$, with $i_k >0$.
It is called {\bf a quasi-line} if all
$i_k=1$.

\hfill

\claim\label{_rationally_conne_Claim_}
Let $M$ be a compact complex manifold containing a 
an ample rational line. Then any $N$ points $z_1, ..., z_N$ can
be connected by an ample rational curve.

\hfill

{\bf Proof:} This fact is well known in algebraic geometry
(see \cite{_Kollar:curves_}). However, its proof is valid
for all complex manifolds. \endproof

\hfill

\claim
Let $M$ be a Riemannian 4-manifold, $\Tw(M)\stackrel \sigma\arrow M$ its twistor
space, $m\in M$ a point, and $S_m:= \sigma^{-1}(m)= S\Lambda^+_m(M)$
the corresponding $S^2$ in $\Tw(M)$. Then $S_m$
  is a quasi-line.

\hfill

{\bf Proof:}
Since the claim is essentially
infinitesimal, it suffices to check it when $M$ is flat.
Then $\Tw(M)= \Tot (\calo(1)^{\oplus 2})
\cong \C P^{3} \backslash \C P^{1}$, and
$S_m$ is a section of $\calo(1)^{\oplus 2}$.
\endproof

\hfill

\corollary 
Any $N$ points $z_1, ..., z_N$ on a twistor space can
be connected by an smooth, ample rational curve
\endproof

\subsection{Rational curves and plurinegative metrics}

For other applications of Gromov's compactness theorem on
manifolds with pluriclosed metrics, 
please see \cite{_Ivashkovich_}.

\hfill

\definition
Let $S$ be a complex curve on a Hermitian
manifold $(M,I,g,\omega)$. Define {\bf  the
Riemannian volume} as $\Vol(S):=\int_S \omega$.

\hfill

\definition
A Hermitian form $\omega$ is called {\bf plurinegative}
({\bf pluripositive})
if the (2,2)-form $dd^c \omega$ is negative (positive).

\hfill

\remark
The notion of ``positive $(k,k)$-form'' comes in two flavours:
{\bf weakly positive} and {\bf strongly positive}
(\cite{_Demailly:ecole_}). When $k=1$ or $k=\dim_\C M-1$,
these two notions coincide. Since in the present paper we
are interesned mostly in 3-dimensional complex manifolds,
this distinction becomes irrelevant. For the sake of a definition,
we shall consider, in the present paper, ``positive'' as
a synonym to ``strongly positive''.

\hfill

Of course, pluriclosed Hermitian
metrics are both pluripositive
and plurinegative.

\hfill

As shown in \cite[(8.12)]{_NHYM_}, 
a standard Hermitian form on a twistor space
of a hyperk\"ahler manifold is pluripositive.

\hfill

\claim
Let $X$ be a component of the moduli of complex curves
on a given  complex manifold, $\tilde X$ the set of pairs
$\{ S\in X, z\in S\subset M\}$, (``the universal family''), and 
$\pi_M:\; \tilde X \arrow M$, $\pi_X:\; \tilde X \arrow X$
the forgetful maps. Then the volume function
$\Vol:\; X\arrow \R^{>0}$ can be expressed as
$\Vol = (\pi_X)_*\pi_M^* \omega$.
\endproof

\hfill

\remark
Since
pullback and pushforward of differential forms
commute with $d$, $d^c$, this gives $dd^c
\Vol=(\pi_X)_*\pi_M^* (dd^c\omega)$.
Therefore, $-\Vol$ is plurisubharmonic on $X$
whenever $\omega$ is plurinegative.

\hfill

\theorem \label{_Gromov_compactness_Theorem_}
(Gromov)
Let $M$ be a compact Hermitian almost complex manifold,
${\goth X}$ the space of all complex curves on $M$, and
${\goth X}\stackrel \Vol \arrow \R^{>0}$ the volume
function. Then $\Vol$ is proper (that is, preimage of
a compact set is compact).

{\bf Proof:}  \cite{_Gromov:curves_}, \cite{_Audin_}.
\endproof

\hfill

\corollary\label{_compa_from_plurine_Corollary_}
Let $M$ be a complex manifold, equipped with a plurinegative
Hermitian form $\omega$, and $X$ a component of the moduli of
complex curves. Then the function $\Vol:\; X \arrow \R^{>0}$
is constant, and $X$ is compact.

\hfill

{\bf Proof:}
Since $\Vol\geq 0$, the set
$\Vol^{-1}(]-\infty, C])$ is compact for all $C\in \R$,
hence $-\Vol$ has a maximum somewhere in $X$. However,  a
plurisubharmonic function which has a maximum is necessarily
constant by E. Hopf's strong maximum principle. Therefore,
$\Vol$ is constant: $\Vol=A$. Now,
compactness of $X= \Vol^{-1}(A)$ follows from Gromov's theorem. \endproof

\subsection{An aside: rational lines on the twistor space of
a K3 surface}

For a complex manifold $Z$ equipped with a 
pluripositive Hermitian form, the same argument
implies that any component of the moduli of 
curves on $Z$ is pseudoconvex.
In particular, this is true on twistor
spaces of hyperk\"ahler manifolds
(\cite{_NHYM_}, \cite{_DD_Mourougane_}).
For a twistor space of K3, a stronger result can be achieved. 

\hfill

\theorem\label{_twi_K3_Stein_Theorem_}
Let $M$ be a K3 surface equipped with a hyperk\"ahler
metric, and $\Tw(M)$ its twistor space. Denote
by $X$ a connected component of the moduli of
rational curves on $\Tw(M)$. Then $X$ is Stein.

\hfill

{\bf Proof:} The proof is based on the following
useful theorem of Forn\ae ss-Narasimhan.

\hfill

\definition
Let $X$ be a complex variety (possibly singular),
and $\phi:\; X \arrow [-\infty, \infty[$ an upper
semicontinuous function. We say that $\phi$ 
is {\bf plurisubharmonic (in the weak sense)} if for any 
holomorphic map $D \stackrel f \arrow X$ from a disc in $\C$,
the composition $f\circ \phi:\; D \arrow \R$ is
plurisubharmonic (or identically $-\infty$). This function is called
{\bf strongly plurisubharmonic} if any petrurbation of
$\phi$ which is small in $C^2$-topology  remains
plurisubharmonic.\footnote{To define precisely what it
means ``small in $C^2$-topology'', we embed 
an open subset $U\subset X$ 
to $\C^n$. Suppose that there exists $\epsilon>0$ such that
for any function $f$ on $\C^n$ with $|f|_{C^2}<\epsilon$,
the sum of $\phi+ f\restrict U$ is 
plurisubharmonic. Then  $\phi$ is called
{\bf strongly plurisubharmonic in $U$}. If $X$ admits
a covering by such $U$, then $\phi$ is called 
{\bf strongly plurisubharmonic in $X$}.}

\hfill

\theorem
Let $X$ be a complex variety admitting an exhaustion function
which is strictly plurisubharmonic. Then $X$ is Stein.

{\bf Proof:} \cite[Theorem 6.1]{_Fornaess_Narasimhan_}.
\endproof

\hfill

Now we can prove \ref{_twi_K3_Stein_Theorem_}.
By Gromov's compactness (\ref{_Gromov_compactness_Theorem_}),
$\Vol:\; X \arrow \R$ is exhaustion, and by 
\cite{_NHYM_} it is plurisubharmonic.
It remains to show that this function is 
strictly plurisubharmonic.

Let $\omega$ be the standard Hermitian form on
$\Tw(M)=M \times \C P^1$.
Denote by $\omega_{\C P^1}$ the Fubini-Study
2-form on $\C P^1$, and let $\pi^*\omega_{\C P^1}$
be its lifting to $\Tw(M)$, where
$\pi:\; \Tw(M)\arrow \C P^1$ is the twistor projection.
Then $dd^c \omega= \omega\wedge \pi^*\omega_{\C P^1}$
(\cite[(8.12)]{_NHYM_}).

Let $v\in Z_xX$ be a vector from the Zariski 
tangent cone of $X$. The strict plurisubharmonicity
of $\Vol$ would follow if the second derivative
$\1 \Lie_v\Lie_{\bar v}\Vol$ is positive for all
$v\neq 0$.

Now, let $x=[S]\in X$ be a point represented by a curve
$S$. Then $Z_x X$ is a subspace of $H^0(NS)$, where $NS$ is
the normal sheaf of $S$. A priori, $S$ can have several
irreducible components, some of them sitting in the
fibers of the twistor projection $\pi:\; \Tw(M)\arrow \C P^1$, 
others transversal to these fibers. However, all
components sitting in the fibers of $\pi$ are fixed,
because the rational curves on K3 are fixed.
Therefore, $v$ in non-trivial along the fibers of $\pi$.
Now, 
\[ \1 \Lie_v\Lie_{\bar v}\Vol=dd^c\Vol(v,\bar v) = 
   \int_S (dd^c \omega)(v,\bar v) 
   \geq \int_S \pi^*\omega_{\C P^1}\cdot \omega(v,\bar v).
\]
The last integral is positive, because $v$ vanishes on 
those components of $S$ which belong to the fibers of $\pi$, 
hence $v\neq 0$ on a component $S_1$ which is transversal to
$\pi$. Then, 
\[ 
  \int_S \pi^*\omega_{\C P^1}\cdot \omega(v,\bar v)\geq 
  \int_{S_1} \pi^*\omega_{\C P^1}\cdot \omega(v,\bar v),
\]
but this integral is positive, because $\pi^*\omega_{\C P^1}$
is positive on each transversal component of $S$.
This proves that $\1 \Lie_v\Lie_{\bar v}\Vol>0$, implying
strict plurisubharmonicity of $\Vol$
\endproof

\hfill

\remark
The variety $X$, which is shown to be Stein in 
\ref{_twi_K3_Stein_Theorem_}, could be singular
(a complex variety is {\bf Stein} if 
it admits a closed holomorphic embedding to $\C^n$).
However, there is not a single known example
of a singular point in any component
 of the space ${\cal S}(M)$  of
rational curves on $\Tw(M)$, when $M$ is a K3.
It is not hard to see that ${\cal S}(M)$ is smooth
when $M$ is a compact torus. It is not entirely
impossible that it is also smooth for a K3.

\subsection{Quasilines and Moishezon manifolds}

Let $M$ be a compact complex manifold,
and $S\subset M$ an ample rational curve.
Assume that the space of deformations of $S$ in $M$
is compact. From \cite[Theorem 3]{_Campana:rational_lines_} 
it follows that $M$ is Moishezon 
(\cite[Remark 3.2 and Theorem 4.5]{_Campana:twistors_}). 
For the convenience of the
reader, I will give an
independent proof of this result
here.

\hfill

Recall that {\bf a quasiline} is a smooth rational 
curve $S\subset M$ such that its normal bundle
is isomorphic to $\calo(1)^n$. A neighbourhood of
a quasiline shares many properties with a neighbourhood
of a line in $\C P^n$. Heuristically, this can be stated 
as follows.

\hfill

{\bf An imprecise statement.}
Let $S\subset M$ be a quasi-line. Then, for an appropriate
tubular neighbourhood $U\subset M$ of $S$, ``for every two points
$x, y\in U$ close to $S$ and far from each other,  there is a 
 unique deformation of $S$ containing $X$ and $Y$.''

\hfill

More precisely:

\hfill

\claim\label{_quasilines_local_defo_Claim_}
Let $S\subset M$ be a quasi-line. Then, for any sufficiently small
tubular neighbourhood $U\subset M$ of $S$, there exists
a smaller tubular neighbourhood $W\subset U$, satisfying the
following condition. Let $\Delta_S$ be
the image of the diagonal embedding $\Delta_S:\; S \arrow W \times W$.
Then there exists an open neighbourhood $V$ of
$\Delta_S$, properly contained in $W\times W$, such that for any pair
 $(x,y) \in W\times W \backslash V$, there exists a unique deformation
$S'\subset U$ of $S$ containing $x$ and $y$.
\endproof

\hfill

\claim\label{_quasilines_local_jet_Claim_}
A small deformation $S'\subset U$ of $S$ passing through $z\in S$
is uniquely determined by a 1-jet of $S'$ at $z$.
\endproof

\hfill

Both of these claims follow from a general results of
deformation theory: the first cohomology of a normal
bundle $NS$ vanish, hence there are no obstructions to
a deformation, and the deformation space is locally
modeled on a space of sections of $NS$. However, 
since $NS=\calo(1)^n$, any section is uniquely determined
by its values in two different points, or by its 1-jet
at any given point.

\hfill

Further on, we shall need the following simple lemma.

\hfill

\lemma\label{_finite_moishezo_Lemma_}
Let $X \arrow Y$ be a domimant map of complex
varieties, which is finite at a general point.
Assume that $X$ is Moishezon. Then $Y$ is
also Moishezon.

\hfill

{\bf Proof:} Replacing $X$ by its ramified covering, we
may assume that outside of its singularities,
the map $X \arrow Y$ is a Galois covering,
with the Galois group $G$. Then $X/G\arrow Y$
is bimeromorphic. Replacing $X$ by its resolution,
we can also assume that $X$ is projective. Then
$X/G$ is also projective, by Noether's theorem
on invariant rings. \endproof

\hfill

Now we can prove the main result of this subsection
(see also \cite[Theorem 3]{_Campana:rational_lines_}).

\hfill

\theorem\label{_Moishezon_if_quasilines_Theorem_}
Let $M$ be a complex manifold, $S\subset M$ a quasi-line,
and $W$ its deformation space. Assume that $W$ is compact.
Then $M$ is Moishezon.

\hfill

{\bf Proof. Step 1:} 
Let $z\in M$ a point, containing a
quasi-line $S \in W$, $W_z$
the set of all curves $S_1 \in W$ containing $z$,
and $\tilde W_z$ -- the set of all pairs $\{x\in S_1, S_1 \in W_z\}$.
From \ref{_quasilines_local_defo_Claim_},  it follows that the map
$\tilde W_z\arrow M$, $(S,x) \arrow x$ is surjective and
finite at generic point.

{\bf Step 2:} By \ref{_finite_moishezo_Lemma_}, 
it would suffice to prove 
that $\tilde W_z$ is Moishezon.

{\bf Step 3:} After an appropriate
bimeromorphic transform, we may assume that $\tilde W_z \arrow W_z$ is
a smooth, proper map with rational, 1-dimensional fibers.
Then $\tilde W_z$ is Moishezon $\Leftrightarrow$ $W_z$ is
Moishezon. Indeed, the space of cycles in a Moishezon
manifold is Moishezon.

{\bf Step 4:} By \ref{_quasilines_local_jet_Claim_}, the map from $W_z$ 
to ${\Bbb P}T_zM$ mapping a quasi-line to its 1-jet is
also generically finite. Therefore, $W_z$ is Moishezon.
\endproof

\hfill

\corollary\label{_twi_pluriclo_then_Moishezon_Corollary_}
Let $M$ be a twistor space admitting a pluriclosed
(or plurinegative)
Hermitian metric. Then $M$ is Moishezon. \endproof

\section{Pluriclosed and Hermitian symplectic metrics on twistor spaces}

Recall the following classical theorem
of Harvey and Lawson (\cite{_Harvey_Lawson:currents_}).

\hfill

\theorem 
Let $M$ be a compact, complex $n$-manifold.  Then the
following conditions  are equivalent. 
\begin{description}
\item[(i)]
$M$ does not admit a K\"ahler metric.
\item[(ii)] $M$ has a non-zero, positive $(n-1,n-1)$-current
$\Theta$ which is an $(n-1,n-1)$-part of a closed current.
\end{description}
\endproof

\hfill

The same argument, applied to pluriclosed or taming
metrics, brings the following (see also 
\cite{_Alessandrini_Bassanellu:positi_}).

\hfill

\theorem \label{_pluriclo_Herm_symple_Hahn_Banach_Theorem_}
Let $M$ be a compact, complex $n$-manifold.  Then
\begin{description}
\item[(a)] $M$ admits no Hermitian symplectic metrics $\Leftrightarrow$ 
$M$ admits a positive, exact, non-zero $(n-1,n-1)$-current.
\item[(b)] $M$ admits no pluriclosed metrics $\Leftrightarrow$ 
$M$ admits a positive, non-zero, $dd^c$-exact $(n-1,n-1)$-current.
\end{description}

{\bf Proof of (a):} Let $A\subset \Lambda^2(M)$ be a cone
of real 2-forms $\eta$ such that the $(1,1)$-part
$\eta^{1,1}$ is strictly positive. Then $A \cap \ker d=0$
is equivalent, by Hahn-Banach theorem, to existence
of a $(2n-2)$-current vanishing on $\ker d$ (hence, exact)
and positive on $A$, hence of type $(n-1,n-1)$ and positive.

\hfill

{\bf Proof of (b):}
Let $A$ be the same as above. Then $A \cap \ker dd^c=0$
is equivalent, by Hahn-Banach, to existence of a 
$(2n-2)$-current $\Theta$ positive on $A$ (hence, of type
 $(n-1,n-1)$ and positive) and vanishing on $\ker dd^c$. 
A positive $dd^c$-exact $(n-1,n-1)$-current 
clearly vanishes on $\ker dd^c$. It remains
to show, conversely, that $\Theta$ is $dd^c$-exact
whenever $\Theta$ vanishes on $\ker dd^c$. 

Since 
$\ker dd^c$ contains the space $\ker d$, the
current $\Theta$  is exact. Let 
\[ H^{n-1,n-1}_{BC}(M,\R):=\frac{\ker d\restrict
   {\Lambda^{n-1,n-1}(M)}}{dd^c(\Lambda^{n-2,n-2}(M))}
\]
be the {\bf Bott-Chern cohomology group}, and
 \[ H^{1,1}_{AE}(M,\R):=\frac{\ker dd^c\restrict
   {\Lambda^{1,1}(M)}}{\im d\cap \Lambda^{1,1}(M)}
\]
be the Aepli cohomology 
(\cite{_Aepli_}, \cite{_Schweitzer_}).
The exterior multiplication induces a pairing
between these two groups, and it is not hard to see
that they are dual. Since $\Theta$ vanishes on 
$\ker dd^c$, the pairing with its Bott-Chern cohomology class
$\langle[\Theta], \cdot\rangle:\; H^{1,1}_{AE}(M,\R)\arrow \R$
vanishes. Therefore, the class of $\Theta$ in 
$H^{n-1,n-1}_{BC}(M,\R)$ is equal zero, hence
$\eta \in \im dd^c$. \ref{_pluriclo_Herm_symple_Hahn_Banach_Theorem_}
is proved.
\endproof

\hfill

This leads to the following useful proposition.

\hfill

\proposition\label{_twi_pluriclo_then_tamed_Proposition_}
Any twistor space $M$ which admits a pluriclosed metric
also admits a Hermitian symplectic structure.

\hfill

{\bf Proof:} 
By \ref{_twi_pluriclo_then_Moishezon_Corollary_},
$M$ is Moishezon. Then, \cite{_DGMS:Formality_}
implies that $M$ satisfies $dd^c$-lemma. Therefore,
any exact (2,2)-current is $dd^c$-exact.
Applying \ref{_pluriclo_Herm_symple_Hahn_Banach_Theorem_},
we obtain that $M$ is Hermitian symplectic.
\endproof

\hfill

\corollary\label{_plurica_then_Kah_Corollary_}
Let $M$ be a twistor space admitting a pluriclosed (or
Hermitian symplectic) metric. Then $M$ is K\"ahler.

\hfill

{\bf Proof:} Th. Peternell (\cite[Corollary 2.3]{_Peternell_}) has 
shown that any compact 
non-K\"ahler Moishezon $n$-manifold admits an exact, 
positive $(n-1,n-1)$-current. Therefore, it is never
Hermitian symplectic (\ref{_pluriclo_Herm_symple_Hahn_Banach_Theorem_}).
By \ref{_twi_pluriclo_then_tamed_Proposition_}, $M$ cannot
be pluriclosed. \endproof

\hfill

{\bf Acknowledgements:}
Many thanks to Gang Tian for  asking
the questions which lead to the writing of this paper, 
and to Frederic Campana and Thomas Peternell
for insightful comments and the reference. Part of this
paper was written during my visit to the
Simons Center for Geometry and Physics; I wish
to thank SCGP for hospitality and stimulating atmosphere.
I am grateful to Nigel Hitchin for the kind interest and 
the inspiration.

\hfill

{\small

\hfill

\noindent {\sc Misha Verbitsky\\
Laboratory of Algebraic Geometry, \\
Faculty of Mathematics, National Research University HSE,\\
7 Vavilova Str. Moscow, Russia
 }

}

\end{document}